\documentclass{amsart}
\usepackage{amsfonts,amscd,amssymb,amsmath,amsthm}
\usepackage{graphicx}

\begin{document}

\newtheorem{theorem}{Theorem}[section]
\newtheorem{lemma}[theorem]{Lemma}
\newtheorem{corollary}[theorem]{Corollary}
\newtheorem{conjecture}[theorem]{Conjecture}
\newtheorem{cor}[theorem]{Corollary}
\newtheorem{proposition}[theorem]{Proposition}
\theoremstyle{definition}
\newtheorem{definition}[theorem]{Definition}
\newtheorem{example}[theorem]{Example}
\newtheorem{claim}[theorem]{Claim}

\newcommand{\R}{\mathbb{R}}
\newcommand{\config}{{\bf C}}

\providecommand{\vect}[1]{\boldsymbol{#1}}
\providecommand{\matr}[1]{\boldsymbol{#1}}

\providecommand{\abv}[1]{\ensuremath{^ \mathrm{#1}}}
\providecommand{\blw}[1]{\ensuremath{_ \mathrm{#1}}}
\providecommand{\unit}[1]{\ensuremath{\mathrm{\, #1}}}

\providecommand{\diff}{\mathrm{d}}
\providecommand{\deriv}[3][ ]{\frac{\diff^{#1}#2}{\diff #3^{#1}}}
\providecommand{\pderiv}[3][ ]{\frac{\partial^{#1}#2}{\partial #3^{#1}}}

\providecommand{\abs}[1]{\lvert#1\rvert}
\providecommand{\norm}[1]{\lVert#1\rVert}

\title{Computational topology for configuration spaces of hard disks}
\author[G. Carlsson]{Gunnar Carlsson}
\address{Dept. of Mathematics, Stanford University}
\email{ <gunnar@math.stanford.edu>}
\author[J. Gorham]{Jackson Gorham}
\address{Dept. of Statistics, Stanford University}
\email{<jacksongorham@gmail.com>}
\author[M. Kahle]{Matthew Kahle}
\address{Dept. of Mathematics, The Ohio State University}
\email{<mkahle@gmail.com>}
\thanks{MK was supported in part by NSA Grant \# H98230-10-1-0227}
\author[J. Mason]{Jeremy Mason}
\address{School of Mathematics, Institute for Advanced Study}
\email{<jkmason@math.ias.edu>}

\date{\today}
\maketitle

\begin{abstract}  We explore the topology of configuration spaces of hard disks experimentally, and show that several changes in the topology can already be observed with a small number of particles.  The results illustrate a theorem of Baryshnikov, Bubenik, and Kahle that critical points correspond to configurations of disks with balanced mechanical stresses, and suggest conjectures about the asymptotic topology as the number of disks tends to infinity.
\end{abstract}

\section{Introduction}
\label{Introduction}

Configuration spaces of hard spheres are of intrinsic interest mathematically and physically, yet so far little seems known about their topology.  As shown in this note, the dependence on the volume fraction (or equivalently on the sphere radius) of properties as fundamental as the connectivity of configuration space is already quite complicated even for a small number of spheres.  More generally, we are concerned not only with the property of connectivity or with the number of connected components $\beta_0$, but also with the $i$th Betti number $\beta_i$ which denotes the dimension of $i$th homology, or colloquially ``the number of $i$-dimensional holes \cite{Hatcher}.''

\begin{figure}
\begin{centering}
\includegraphics[width=3.25in]{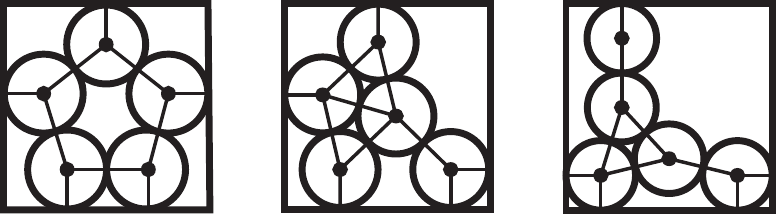}
\bigskip
\includegraphics[width=3.25in]{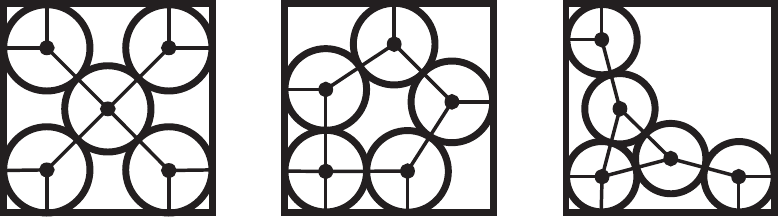}\\
\end{centering}
\caption{Index zero (bottom row) and index one (top row) critical points, shown with their bond networks.}
\label{fig:mins}
\end{figure}

One purpose of this article is to illustrate some of the complexity of the topology of configuration spaces of hard spheres, even in two dimensions (i.e.\ hard disks) when the number of disks is small.  We illustrate a theorem of Baryshnikov, Bubenik, and Kahle that mechanically-balanced (or jammed) configurations act as critical points, indicating the only places where the topology can change \cite{BBK}.  We always mean ``critical point'' in the mathematical sense of non-differentiable or singular, rather than the statistical mechanical sense.  Examples of such configurations may be seen in Figure \ref{fig:mins}.  The combinatorial and geometric structure of the bond networks of these critical points is important in classifying their degeneracy and in computing their index.  Furthermore, the bond networks suggest connections to the geometric structure approach in hard particle packing \cite{Sal}, as well as to tensegrity and strut frameworks \cite{Conn}.

We show here that the topology of these configuration spaces can already be fairly complicated in two dimensions with a small number of particles.  For example, with five disks in the unit square, the topology changes at least $20$ times as the disk radius varies, and for radius $0.1686 < r < 0.1692$ the configuration space has $\beta_1=2176$ ``holes''.

The topology of configuration spaces has been studied in the context of motion planning for robots; see the surveys by Farber and Ghrist \cite{Farber, Ghrist}, and the textbook by LaValle \cite{Lavalle}.  One central idea in this field is that topological complexity, measured in various ways, is a measure of the difficulty of constructing a motion planning algorithm.  Many papers idealize robots as point particles in a suitable metric space --- although it has been noted that the assumption that they have no volume is not physically realistic \cite{AG}.  On the other hand, expanding the particles to have positive thickness greatly complicates the topology of the underlying configuration space --- see for example the recent work of Deeley on thick particles on a metric graph \cite{Deeley}.

%The configuration spaces of $n$ hard disks of radius $r$ in a unit square is already fairly mysterious.  As Persi Diaconis noted,
%\begin{quote}
% ``We know very, very little about the topology of the set of configurations: for fixed $n$, what are useful bounds on $r$ for the space to be connected? What are the Betti numbers? Of course, for $r$ small this set is connected but very little else is known.'' \cite{Persi_Markov}.
%\end{quote}
%It was shown recently

As further motivation for this study, a set of hard disks inhabiting a bounded area may be regarded as a model thermodynamic system.  From this standpoint, a topological property such as connectivity of the configuration space is a fundamental concern, since this is required to satisfy the assumptions of equal a priori probabilities and ergodicity (e.g., Section VII of Torquato and Stillinger \cite{Sal}) inherent to the statistical mechanics formalism.  Despite the importance of this question, few results seem to be known for finite hard disk systems \cite{Persi_Markov}, and fewer still in the thermodynamic limit.  It is known that if one can only move one disk at a time, the configuration space for large numbers of disks can be disconnected even for density arbitrarily close to zero \cite{Boro, Still}.  Even if collective motions of disks are allowed, the configuration space may be disconnected for density bounded away from maximal --- for example, the reinforced Kagome lattice is strictly jammed and only has density $\approx 0.68$ compared to the close packing density of $\approx 0.91$ \cite{Donev}.

Perhaps one of the most notable properties of the hard disk system in the context of thermodynamics is the existence of a solid-to-liquid phase transition \cite{Alder1962} with the variation of the packing fraction, or equivalently of the disk radius.  This is particularly suggestive in the context of a number of recent papers that explore the hypothesis that statistical-mechanical phase transitions are intimately related to changes in the topology of surfaces of constant potential (i.e. equipotential submanifolds) in the configuration space \cite{Teix, Grinza, exact,  FF}.  Specifically, a theorem of Franzosi, Pettini, and Spinelli \cite{FPS1,FPS2} states that for rather benign conditions (smooth, finite-range, confining potentials), the low-order derivatives of the Helmholtz free energy cannot display a discontinuity unless there is simultaneously a change in the topology of the equipotential submanifolds with the variation of a relevant parameter, such as the disk radius.  While this indicates that a topological change is necessary, a conjecture known as the Topological Hypothesis roughly states that a major topological change (a proposal for specific conditions is given by Kastner, Schreiber and Schnetz \cite{Kastner2007, Kastner2008}) is sufficient to drive the phase transition as well.  This immediately raises the question of the topological signature of the solid-to-liquid phase transition for the hard disk system.

Historically, the accumulation of certain indirect evidence encouraged the development of the Topological Hypothesis.  For example, the magnitude of fluctuations in the curvature of equipotential submanifolds was observed to show singular behavior near a phase transition \cite{Caiani1997}.  This singular behavior seemed to be insensitive to the choice of metric \cite{Franzosi1999}, pointing to a more fundamental change in the topology of the accessible portion of configuration space at the phase transition.  Further support was provided by the existence of major changes in the topology of the equipotential submanifolds accompanying phase transitions in certain exactly solvable systems \cite{Casetti1999, Casetti2002, Grinza, exact}.  Strikingly, the presence of these changes is known to persist for small systems well away from the thermodynamic limit \cite{Casetti1999, Franzosi2000}.  That is, a major topology change in the equipotential submanifolds of a system with just a few hard disks may indicate the presence of a topology change that drives the phase transition of the corresponding system in the thermodynamic limit.

There is a persistent difficulty faced by investigations of the Topological Hypothesis in the physics literature, namely, that the exactly solvable systems mentioned above do not satisfy the conditions for the theorem of Franzosi, Pettini and Spinelli.  Specifically, they violate the requirements of short-range interactions \cite{Casetti1999, Casetti2002, exact} and confining potentials \cite{Grinza}.  This stands in sharp contrast to the short-range, hard-core repulsive potential characteristic of hard disks systems, meaning that a complete characterization of the configuration space topology of our system will begin to address a need that is openly acknowledged in the literature \cite{exact}.

%In other words, for many systems a change in topology is necessary for a phase transition.  However even for these systems the converse fails, and topological changes are clearly not sufficient to guarantee phase transitions.  Indeed, one can give examples of systems where the number of topological changes is comparable to the number of degrees of freedom, so is tending to infinity with the number of particles, but where there is only one statistical mechanical phase transition.

%On the other hand there are also many interesting systems for which Franzosi, Pettini, and Spinelli's theorem does not apply, and examples have been given by Risau-Gusman, Ribeiro-Teixeira, and Stariolo \cite{short-range}, and by Kastner \cite{Kastner}, of models where topology alone can not capture the phase transitions. However in these systems, the topology is fairly uncomplicated (e.g.\ contractible manifolds).  One might hope that the Topological Hypothesis holds in systems where the topology of configuration spaces is sufficiently complicated or changes a large number of times.

We begin this study by reconstructing the dendrogram, i.e.\ a graph which encodes the appearance of path components and the connections between them as the disk radius (or volume fraction) varies.  More generally, for a suitably defined energy function, the dendrogram is defined to be the graph with nodes given by the basins of attraction in the energy landscape and edges given by the saddle points that connect the basins as the energy increases.  These correspond in the physics literature to collectively-jammed configurations and to the transitions between them (see Torquato \& Stillinger \cite{Sal}), subjects of continuing interest to the physics community.  The name ``dendrogram'' comes from hierarchical clustering in statistics.  We find the saddle points using the ``nudged elastic band'' method \cite{h2000}.  From the point of view of Morse theory, the index of the basins of attraction is zero since there is no collective motion that decreases the energy, and the index of the saddle points is one since there is one dimension of collective motion that decreases the energy.  Examples of these varieties of critical points appear in Figure \ref{fig:mins}.

As alluded to above, these features of the configuration space correspond in the Morse-theory literature to critical points of index zero and index one, respectively \cite{Milnor}.  Along with the notion of a critical index comes the notion of higher index critical points, allowing for the construction of a more complete picture of the configuration space and the related energy landscape.  This is precisely the field of Morse theory, which provides a framework to pass from the populations of critical points with specified indices to a reconstruction of the configuration space, and allows us to answer many more questions than simply the energy at which the saddle point connecting jammed configurations is accessible.   For instance, we deduce certain facts about the Betti numbers $\beta_k$ that indicate the number of $k$-dimensional holes or the dimension of $k$th homology, and observe suggestive patterns about the asymptotic topology of configuration spaces for large numbers of disks.

There is good reason to introduce Morse theory as a fundamental instrument for the study of configuration spaces.  Various flavors of Morse theory sit at the heart of Franzosi, Pettini, and Spinelli's proof \cite{FPS1} in the context of the Topological Hypothesis, Yao et al.'s  exploration of low-density states in biomolecular folding pathways \cite{Yuan}, and Deeley's study of configuration spaces of thick particles on a metric graph \cite{Deeley}.  Section \ref{sec:Morse} will briefly overview Morse theory, though most of the technical mathematical details, and particularly the proof of the main theorem we quote, are left to another article \cite{BBK}.  Our primary focus in this article is on our numerical experiments to identify critical points of various indices and connecting the Morse theoretic viewpoint of configurations spaces with related concepts in the physics literature.

%We study these configuration spaces using a computational implementation of Morse theory.  We note that Morse theory is also at the heart of Franzosi, Pettini, and Spinelli's proof. Our main concern in this article is to understand the structure of connected components, and how this structure changes as the radius of disks varies.  Balanced configurations play the role of critical points in Morse theory.

\section{Configuration spaces}

Let $\config(n,r)$ denote the configuration space of $n$ nonoverlapping disks of radius $r$ in the unit square $[0,1]^2$.  If we indicate the coordinates of the center of the $i$th disk by $x_i$, we can write this as
$$ \config(n,r) = \{(x_1, x_2, \dots, x_n) \mid x_i \in [r, 1-r]^2, d(x_i, x_j) \ge 2r \}.$$
Since $\config(n,r) \subset \R^{2n}$, it naturally inherits a subspace topology and a metric from the Euclidean space of the appropriate dimension.  The configuration space of hard disks in the limit of large $r$ is the empty space, since the disks do not fit in the unit square.  Meanwhile, the configuration space of $n$ labelled points in the plane (or more generally on a manifold) is well studied \cite{Cohen,Cohen2}, and is equivalent to the configuration space of hard disks in the limit of small $r$.  This leaves the study of the topology of $\config(n,r)$ for intermediate $r$, which is much more complicated than for the two cases above.

Notice that if $r' < r$, then there exists a natural inclusion map $$\config(n,r) \hookrightarrow \config(n,r')$$ given by shrinking the disks with fixed centers.  In other words, as the disk radius shrinks, the space of allowable configurations grows.

As mentioned in Section \ref{Introduction}, a basic question is whether $\config(n,r)$ is connected for a given value of $r$, and if not, how many components there are.  A more subtle description involves the construction of the dendrogram, or the graph that encodes the appearance of persistent components and the connections between them as the disk radius varies.  Equivalently, the dendrogram is a means to visualize distinct basins of attraction in the energy landscape and the saddle points between them, not the properties of the configuration space at higher energies.  That would require identifying the critical points of all indices and using Morse theory to reconstruct the remaining features of the configuration space.

We give a brief discussion of Morse theory to establish a unifying set of language and ideas.  Following that, we apply each of these approaches to the specific case of the configuration space of five disks in the unit square.

\section{Morse theory}
\label{sec:Morse}

The basic idea of Morse theory is to filter a topological space $X$ by the sub-level sets $f^{-1}( (-\infty,r]) $ of a suitable function $f: X \to \R$ \cite{Milnor,Forman}.  Suppose that the function $f$ assigns an energy $E$ to every point in the configuration space, and that this energy increases as the radius of the disks decreases.  The sub-level set associated with a particular energy $E_0$ is just the set of accessible states, or the set of points in the configuration space where the energy is less than or equal to $E_0$.  The existence of the inclusion map in the preceding section is equivalent to the statement that as $E_0$ increases, the set of accessible states increases.  Provided that the function $f$ is well-behaved, the topology of the sub-level sets will change only a small number of times as $r$ decreases.  One of the tenets of traditional Morse theory is that these changes occur precisely at the critical points where $\nabla f(x) = 0$, or at the mechanical equilibria of the system.  Using techniques similar to those in ``Min-type'' Morse theory \cite{min-type} and stratified Morse theory \cite{Goresky}, Baryshnikov, Bubenik, and Kahle recently proved the following theorem \cite{BBK}.

\begin{theorem}\label{thm:bbk} If there are no ``balanced'' configurations of disks of radius $r$ with $r \in [a,b]$ then $\config(n , a)$ is homotopy equivalent to $\config(n,b)$.
\end{theorem}

The definition of a ``mechanically-balanced'' configuration is as follows. For a configuration of disks, define the {\it bond network} to have vertices at disk centers, and edges $e_{ij}$ connecting pairs of tangent disks $\{i,j \}$, and also disks to the boundary of the region.  We consider the edge $e_{ij}$ to be a vector, pointing from $i$ to $j$, so in particular $e_{ij} = - e_{ji}$.

\begin{definition} A configuration is {\it mechanically-balanced} if non-negative weights $w_{ij}$ can be assigned to the edges of the bond network so that at least one $w_{ij}$ is positive, and 
for every fixed vertex $i$, $$\sum_{j \sim i} w_{ij} e_{ij} =0 .$$
where the sum is over all vertices $j$ adjacent to $i$.
\label{balanced}
\end{definition}

The point here is that $\config(n,r)$ is the sub-level sets of some function --- but not a smooth function.  Smooth Morse theory and Morse-Bott theory don't quite work in this case, but something similar to min-type Morse theory \cite{min-type} does work.  We note that this is similar in spirit to the PL Morse-Bott Theory in Deeley's recent work \cite{Deeley}.

For a smooth Morse function the index of a non-degenerate critical point is the number of negative eigenvalues of the Hessian.   For a min-type Morse function (such as for hard disks), you can compute the index of critical points with linear programming, i.e.~the index is the dimension of a certain cone.  In either form of Morse theory, the fundamental theorem states that every time you have an index $k$ critical point, the homotopy type changes by gluing in another $k$-dimensional cell, and importantly --- that there are not any other changes in the topology except when there are critical points.

\section{Part I:  Dendrogram}

It is convenient for our methods to first slightly ``soften'' the disks.  We recover their hardness shortly.   
In general, on the configuration space of $n$ distinct points in a region $R$ we can define a smooth energy function $E$ by
$$E =\sum_{1 \le i < j \le n}  \frac{1}{(d(x_i, x_j)/2)^h}+  \sum_{1 \le i \le n} \frac{1}{d(x_i, \partial R)^h}  $$
%
%$$E =\sum_{1 \le i < j \le n} \frac{2^h}{\left( (x_i-x_j)^2 + (y_i-y_j )^2  \right) ^{h/2} }+\sum_{1 \le i \le n} \left( \frac{1}{x_i^h}+\frac{1}{(1-x_i)^h}+\frac{1}{y_i^h}+\frac{1}{(1-y_i)^h} \right),$$
where $h \gg 0$ is a fixed ``hardness'' parameter and $\partial R$ denotes the boundary of $R$.  For our experiments we set $h=50$, $n=5$, and $R=[0,1]^2$.
%
%Five disks is already enough to give complicated topological behavior as the radius varies.  The choice of square boundary is somewhat arbitrary, and we expect one would find qualitatively similar results for a different shaped convex region or even for a region without boundary such as a torus or sphere.  However for a region without boundary (or even a smooth boundary such as a disk), one will likely have to add a few more disks before the complications appear.

%\subsection{Dendrogram}

Taking uniform random points in $[0,1]^{10}$ and flowing in the direction of steepest descent along $- \nabla E$, we find the three qualitatively different types of local minima illustrated in the bottom row of Figure \ref{fig:mins}.  We believe this to be a complete list of all ``types'' of local minima (see discussion below) for five disks in a box.

\begin{figure}
\begin{centering}
\includegraphics[width=4in]{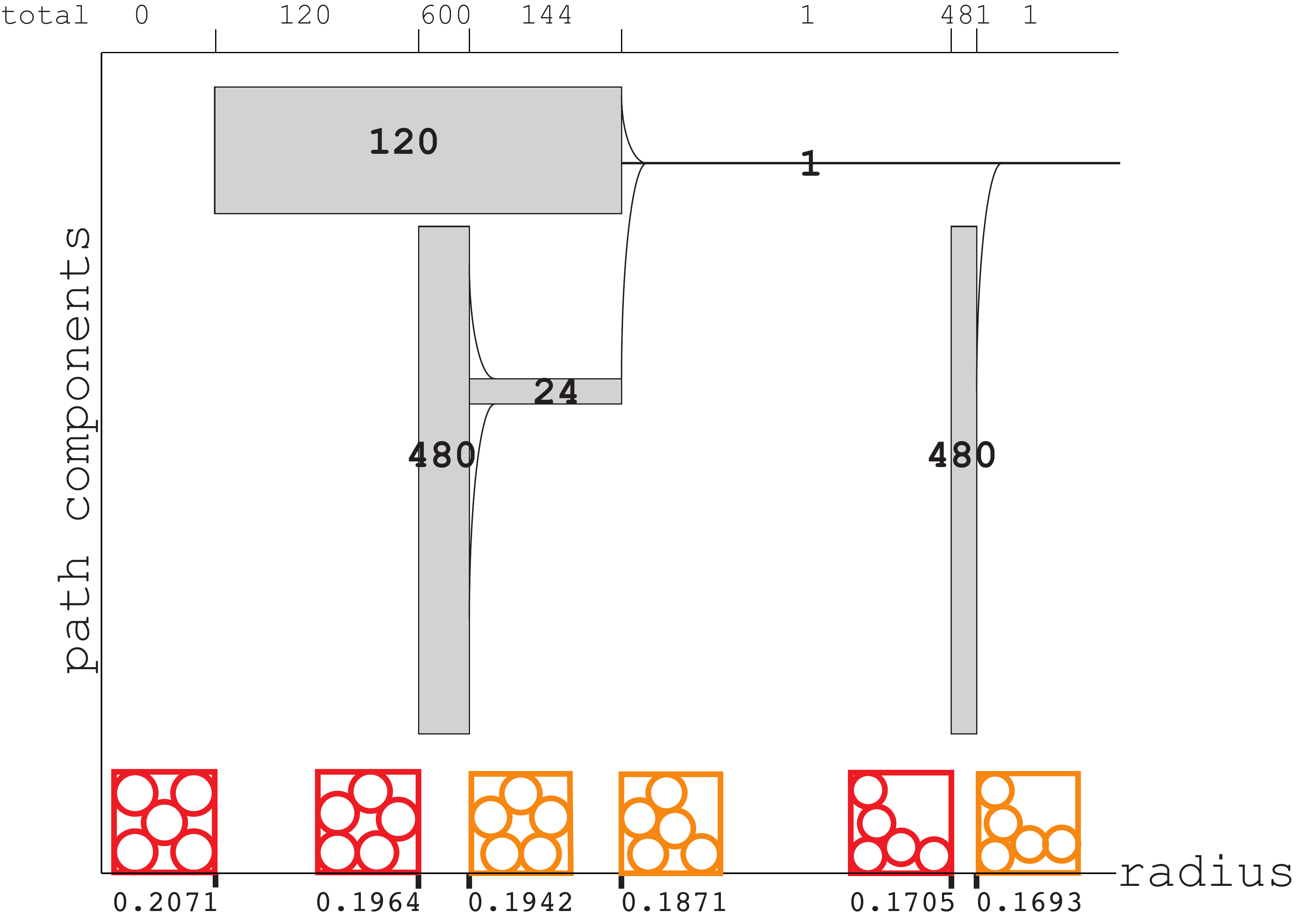}
\end{centering}
\caption{(Color online.)  A bird's-eye view of the dendrogram for labeled disks.  New components appear at radius $0.2071$, $0.1964$, and $0.1705$, and components merge together at radius $0.1942$, $0.1871$, and $0.1693$.}
\label{fig:Reeb1}
\end{figure}

Notice that a point in $[0,1]^{10}$ gives the coordinates of five distinguishable disks.  Since exchanging the labels of the disks does not change the energy of the configuration, this permutation symmetry entails that every critical point be symmetrically equivalent to $120$ points in the configuration space (including itself).  (To mathematicians, the labeled configuration space admits a free action of the symmetric group, and the quotient is the unlabeled configuration space.)  This number of instances of each type is further increased by the additional configurations generated by applying the symmetries of the unit square.  For instance, the class of configurations symmetrically equivalent to the critical point at a radius of $0.2071$ in Figure \ref{fig:Reeb1} contains $120$ points, while the equivalence class containing the critical point at a radius of $0.1964$ contains $480$ points, and the equivalence class containing the critical point at a radius of $0.1693$ contains $960$ points.  So by a ``type'' of critical point, we mean an equivalence class of critical points under both types of symmetries.

Once a local minimum of $E$ is found, the preliminary radius $r$ of the configuration is the maximum possible radius such that that no disks overlap or extend outside the unit square.  The ``bond network'' of the configuration is found by connecting centers of disks that are within $(2 + \epsilon) r$ of each other, or within $(1+\epsilon) r$ of the boundary.  Here $\epsilon > 0$ is a fixed small tolerance parameter, usually around $0.001$.  With the bond network established, we effectively harden the disks by searching for an exact solution to the set of equations given by varying $r$ and setting $\epsilon = 0$.  In some cases this gives a unique solution, and establishes the final values of the radius for the three red configurations along the bottom of Figure \ref{fig:Reeb1}.

\begin{figure}
\begin{centering}
\includegraphics[width=2.5in]{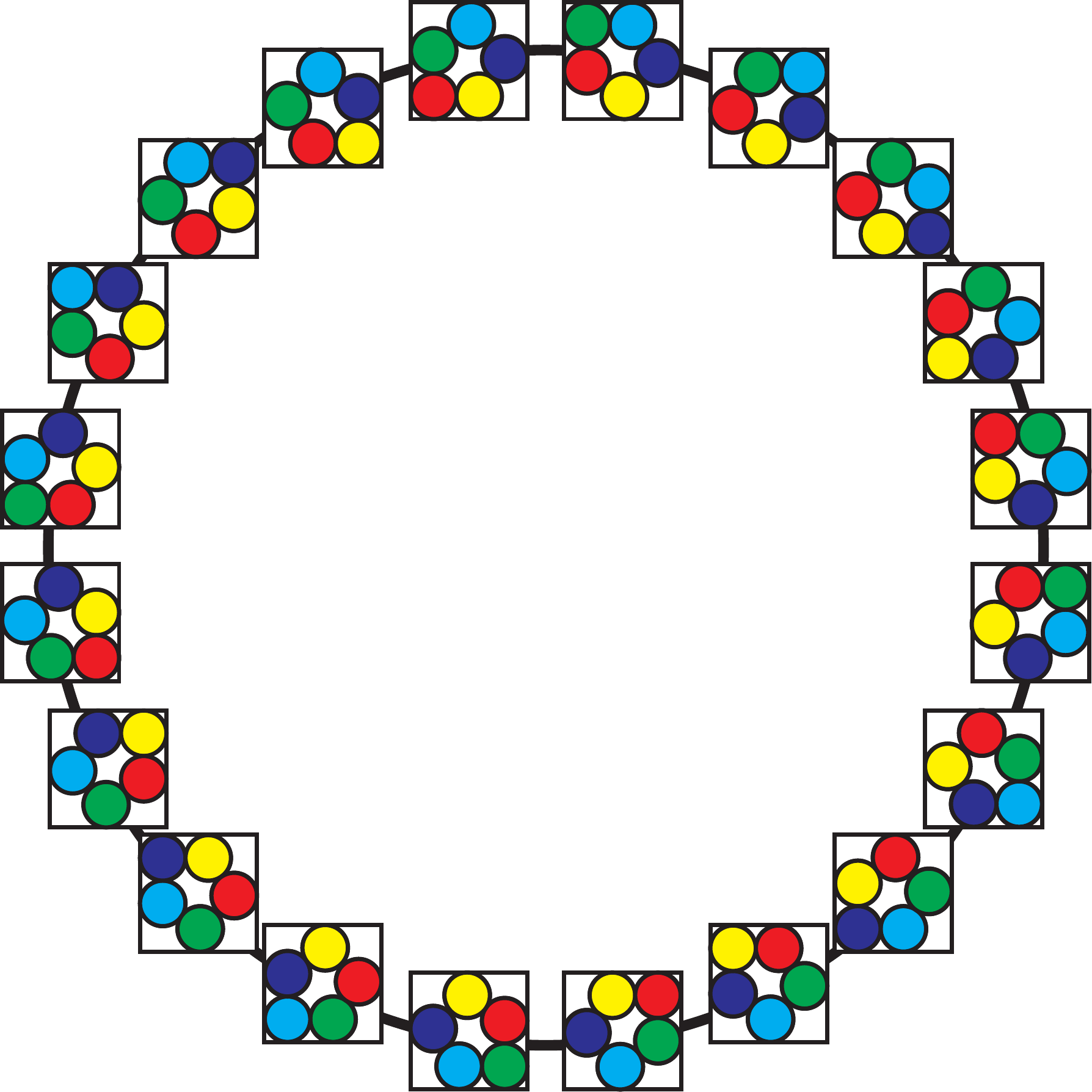}
\end{centering}
\caption{(Color online.)  For radius $ 0.1871 < r < 0.1942$, each connected component of the configuration space is homotopy equivalent to a circle.  One such embedded circle and $20$ index-$0$ critical points is shown here.  Not shown are the $20$ transitionary index-$1$ critical points between them. }
\label{fig:cycle}  
\end{figure}

\begin{figure}
\begin{centering}
\includegraphics[width=2.5in]{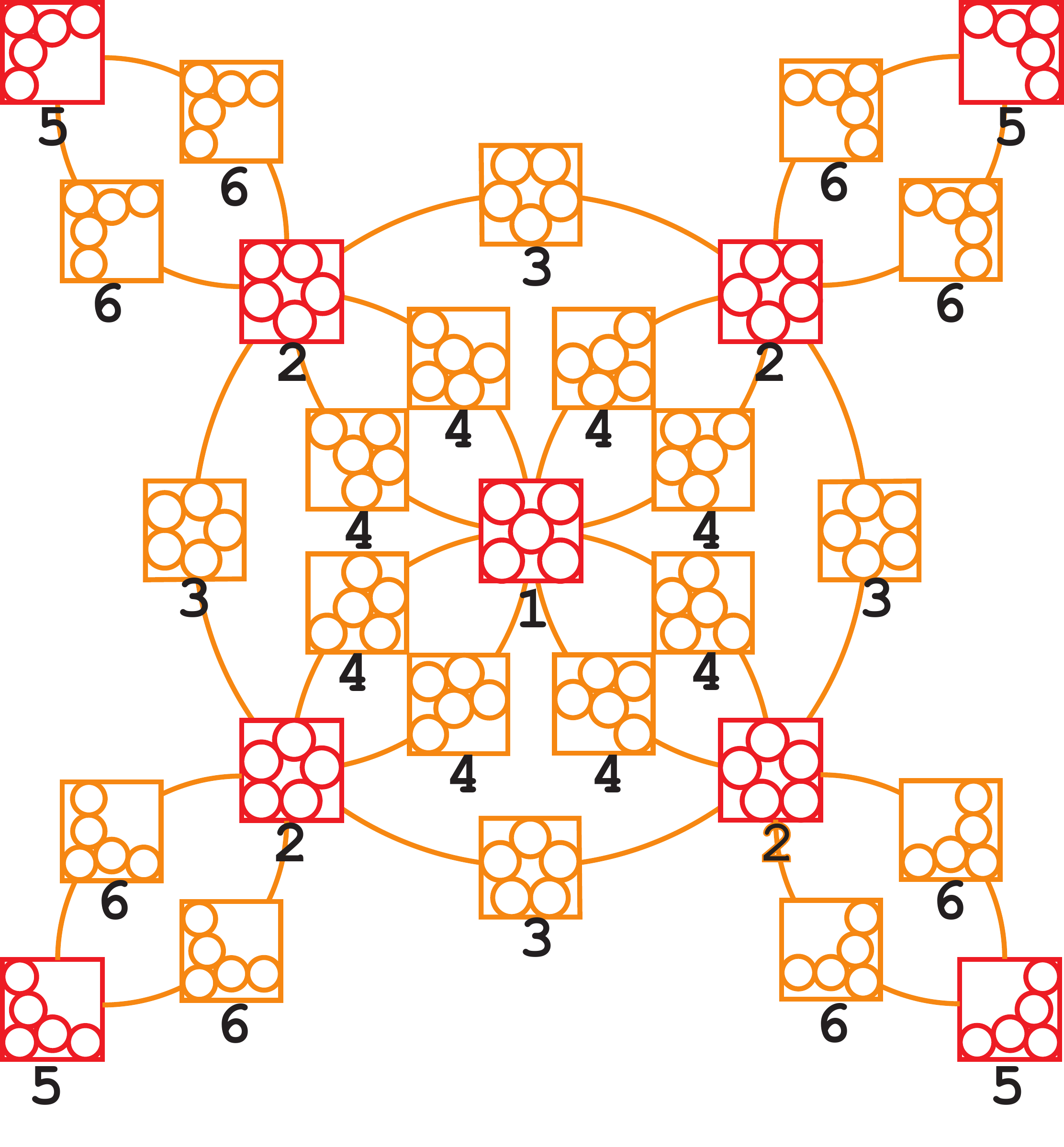}
\end{centering}
\caption{(Color online.) The dendrogram for unlabeled disks.  The numbers indicate the order in which critical points appear as the disks shrink --- new components appear at $1$, $2$, and $5$, and components merge at $3$, $4$, and $6$.  The radius where each vertex and edge appear is given in Figure \ref{fig:index5}. }
\label{fig:Reeb2}
\end{figure}

We then consider pairs of local minima which are close in the configuration space for the induced metric and find low energy paths between them using the ``nudged elastic band'' method \cite{NEB,h2000}.  The highest energy points along a low energy path should be close to a saddle point.  Indeed, by establishing the bond networks for all of the highest energy points along low energy paths and applying the hardening procedure described above, we identify the saddle points depicted in orange along the bottom of Figure \ref{fig:Reeb1}.  This allows us to construct the dendrogram for our system, though a full graphical representation of this graph would be quite complicated.  We present simplified versions of some of the salient features in Figures \ref{fig:Reeb1}, \ref{fig:cycle} and \ref{fig:Reeb2}.

\begin{figure}
\begin{centering}
\includegraphics[width=2in]{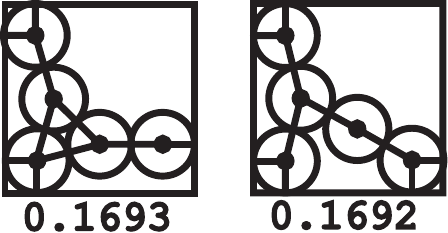}
\end{centering}
\caption{Two saddle points with nearly the same radius}
\label{fig:closer}
\end{figure}

A precise value of the radius is important to be able to order the critical points by the radius of the disks when they appear.  Evidence that a rough numerical approach is not sufficient is provided by the two critical points of index one in Figure \ref{fig:closer}, where the radius differs only in the fourth decimal place.  This is one reason why the more elaborate approach described above is necessary.

%Once we find one local minimum we take advantage of the symmetry of the situation -- by considering all permutations of disks, and also all symmetries of the square, once we find one local minimum we find many.

%When $r = 0.2071$ the first connected components appear, and there are $120$ of them by reordering the disks.  At $r \approx 5.09$ another $480$ components appear.  At $1/r \approx 5.15$ these $480$ components collapse down to $24$.  A close up is shown in Figure \ref{fig:cycle}--- here the disks are colored to emphasize that they are labeled disks, and one sees $20$ jammed configurations arranged along a circle.
%At $1/r \approx 5.34$ the configuration space becomes connected.   When $1/r \approx 5.87$ and $1/r \approx 5.91$, $480$ more components appear and then join the rest of the configuration space.

Figure \ref{fig:Reeb1} shows the variation in the number of connected components as a function of the radius.  The total number of connected components for a given radius is shown at the top of the figure.  The central part of the figure collects the connected components into classes of symmetrically equivalent elements, and labels the classes by the number of elements they contain.  The bottom of the figure shows the configurations of the disks at the relevant critical points, and provides the values of the radius when these configurations appear.

For $r > 0.2071$, the disks do not fit inside the unit square and the configuration space is empty.  The first allowed configurations occur at a radius of $0.2071$, when $120$ connected components appear simultaneously.  All of these belong to the same class, and correspond to different labelings of the disks in the first red configuration.  As the radius decreases, a class of $480$ more connected components appear at radius $0.1964$, corresponding to different labelings and arrangements of the disks in the second red configuration.  The $480$ components of this class connect in groups of $20$ with the emergence of the first saddle point at radius $0.1942$, leaving a class of $24$ embedded circles.  One of these circles is shown in Figure \ref{fig:cycle}.  Every arc segment connecting the jammed configurations in Figure \ref{fig:cycle} is a low energy path containing a saddle point of the type appearing at radius $0.1942$.  The appearance of the second saddle point at radius $0.1871$ connects all of the remaining components.  The space remains connected for all smaller values of the radius, with the exception of the appearance of $480$ more components at a radius of $0.1705$, though these connect up with the main component soon afterward.

An alternate representation of the dendrogram is given in Figure \ref{fig:Reeb2}.  Since the purpose of this figure is to allow the visualization of the paths connecting energy minima, we consider only unlabeled disks for the sake of simplicity.  The labels of the critical points indicate the order of appearance as the radius of the disks decreases, where a set of critical points with the same labels constitutes an equivalence class.  Similar to the description in Figure \ref{fig:Reeb1}, the first allowed configuration to appear is the red critical point in the center of the figure.  As the radius decreases this single configuration is joined by the four critical points arrayed around the central circle, though they all remain isolated until the first saddle points appear.  At this interval the configuration space is given by the isolated configuration in the center and a circle that is roughly analogous to that in Figure \ref{fig:cycle}, but with unlabeled disks.  The central configuration and the circle become connected with the appearance of the second saddle point, and the space remains a single component until the final energy minima appear at the outer extremities of the figure.  The third saddle point once more connects the space into a single component.

Figures \ref{fig:Reeb1} and \ref{fig:Reeb2} do not show all of the critical points of index one for five disks in the unit square.  Indeed, we identified five types of critical points of index one, though only three of these change the number of connected components in the configuration space.  Figures \ref{fig:Reeb1} and \ref{fig:Reeb2} show precisely these three types of critical points, since the remaining two do not directly affect the connectivity of the space.

%\begin{figure}
%\begin{centering}
%\includegraphics[width=3in]{index1n.pdf}
%\end{centering}
%\caption{The relevant saddle points.}
%\label{fig:saddle}
%\end{figure}

\section{Part II: Higher index critical points}

While the above section is concerned with the critical points of index zero and index one and the relation of these to the dendrogram, the current section is concerned with more general features of the configuration space.  As mentioned in Section \ref{Introduction}, this requires the identification of the critical points of higher index.  For this purpose, we work with the same smooth function $E$ defined above.  Using the MuPAD computer algebra system included in the Symbolic Math Toolbox for MATLAB, we symbolically differentiate $E$ and set
\begin{equation*}
F = \norm{\nabla E}^2.
\end{equation*}
Since $\nabla E = 0$ at the critical points of $E$, $F = 0$ at the critical points of $E$ as well and is strictly greater than zero everywhere else.  Hence, we may find critical points of $E$ of any index by flowing down the function $F$ to one of the zeros.

More specifically, our search begins from a random point in $[0, 1]^{10}$.  We symbolically calculate the gradient of $F$ with MuPAD and apply the conjugate gradient method with the Polak--Ribi\`ere formula, using the secant method to find the minimum of $F$ in the search direction on every iteration.  This gives much more rapid convergence to the local minima of $F$ than the steepest descent algorithm.

Once the conjugate gradient algorithm converges, the preliminary radius $r$ of the configuration is the maximum possible radius such that that no disks overlap or extend outside the unit square.  The bond network is constructed just as for the calculation of the dendrogram, though the tolerance parameter $\epsilon$ is on the order of $0.01$.  We recover the configuration of the critical point in the hard disk system by modeling the bonds in the bond network as springs.  Let $s$ be the average of the bond lengths from the disks to the boundary and half the bond lengths connecting pairs of disks.  Set the equilibrium length of the springs from the disks to the boundary to be $s$, and the equilibrium length of the springs connecting pairs of disks to be $2s$.  We minimize the potential energy of the spring network using the conjugate gradient method and update $s$ on every iteration.  This is repeated while reducing the spring constant from a small positive number to zero, allowing us to find the radius of the disks at a critical point of the hard disk system to more than five significant digits.

The index of a critical point is found by symbolically calculating the Hessian of $E$, evaluating this expression for the given disk locations, and numerically finding the number of negative eigenvalues.  Provided that the Hessian is well-conditioned, the index of a critical point calculated for the smooth function $E$ should be the same as for the hard disk system.  Some discretion should be used when finding the critical index this way though, since the Hessian is not always well-conditioned.  In these cases, the index is found as the dimension of a certain cone calculated from the bond network.

\begin{figure}
\begin{centering}
\includegraphics[width=3.5in]{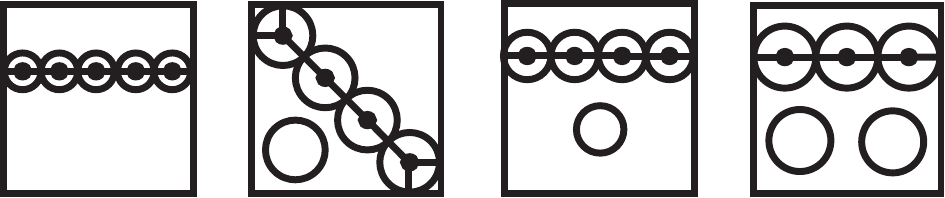}
\end{centering}
\caption{From left to right, critical submanifolds of dimensions $1$, $2$, $3$, and $5$. }
\label{fig:subm}
\end{figure}

We found only $30$ types of critical points and critical submanifolds after repeating this process on the order of $10^5$ times.  By ``critical submanifold'' we mean a critical point with several degrees of freedom, where the number of degrees of freedom is the dimension of the submanifold.  (In the strict mathematical sense, these might not be manifolds or even manifolds with boundary.) Specifically, the presence of each ``rattler'' contributes two dimensions to the critical submanifold, while a row of disks across the square adds one dimension.  Several examples of these appear in Figure \ref{fig:subm}, along with the corresponding dimensions.  Notice that the $5$-dimensional sub-manifold on the right is not contractible, as the others are --- this one is homotopy equivalent to a circle.

\begin{figure}
\begin{centering}
\includegraphics[width=3.5in]{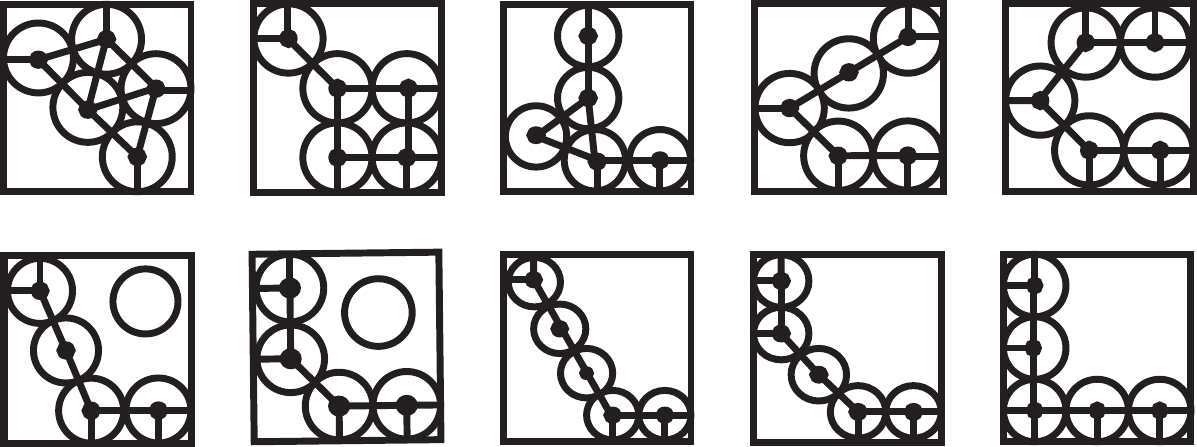}
\end{centering}
\caption{Degenerate types}
\label{fig:degen}
\end{figure}

Of these $30$ types of critical points and submanifolds, we consider the $10$ in Figure \ref{fig:degen} to be degenerate.  Our criterion is if weights $w_{ij}$ may be assigned to the edges of the bond network of a balanced configuration so that at least one $w_{ij}$ is positive and at least one $w_{ij}$ vanishes, then the critical point or submanifold is degenerate.  Computing the local topology near several critical points that satisfy this condition has convinced us that the topology does not change at these critical points, but we do not yet state this as a formal theorem.  We intend to explore the question of degeneracy in more mathematical detail in a future paper.

\begin{figure}
\begin{centering}
\includegraphics[width=4.25in]{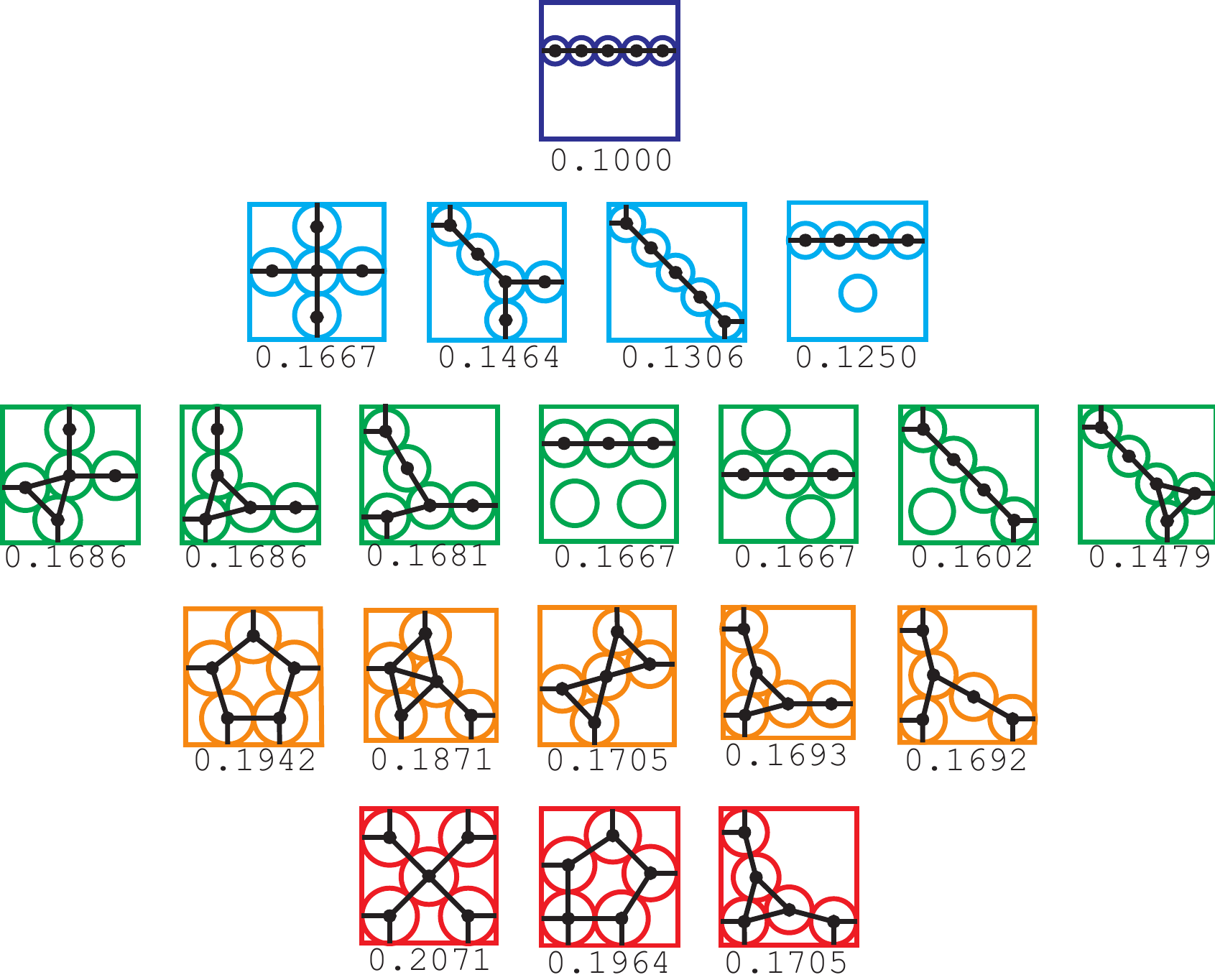}
\end{centering}
\caption{(Color online.) Non-degenerate types, in rows by index.  Index $0$ critical points are in the bottom row, all the way up to index $4$ in the top row.}
\label{fig:index5}
\end{figure}

The remaining $20$ non-degenerate critical points and submanifolds appear in Figure \ref{fig:index5}, ordered horizontally by the radius of the disks and vertically by index.  

This figure reveals that higher index critical points generally coincide with smaller radius.   For instance, all of the index-$0$ and index-$1$ critical points occur at larger radius than any of the higher-index critical points.   The fundamental theorem in Morse theory states that each critical point of $i$ contributes an $i$-dimensional cell via an attaching map, up to homotopy equivalence (smooth deformation) \cite{Goresky, Milnor}.  This indicates that for radius $r \ge 0.1692$, the whole configuration space is homotopy equivalent to a $1$-dimensional cell complex, i.e.\ a graph.

\begin{figure}
        \begin{tabular}{r|rrrrrrrrr}
       && \includegraphics[width=0.3in]{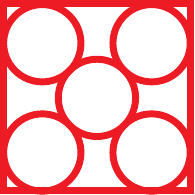}&\includegraphics[width=0.3in]{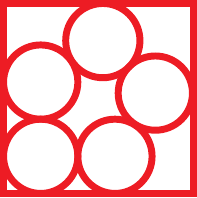}&\includegraphics[width=0.3in]{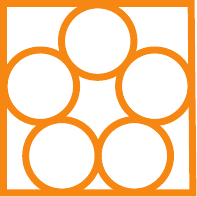}&\includegraphics[width=0.3in]{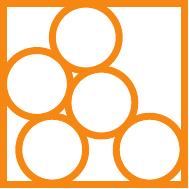}&\includegraphics[width=0.3in]{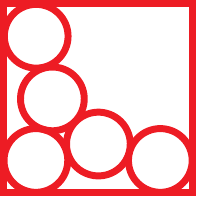}&\includegraphics[width=0.3in]{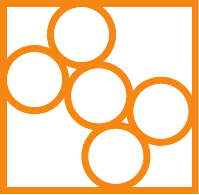}&\includegraphics[width=0.3in]{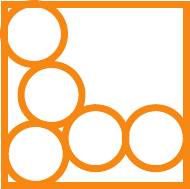}&\includegraphics[width=0.3in]{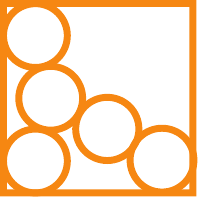}\\
       %\hline
        radius &  & 0.2071 & 0.1964 & 0.1942 & 0.1871 & 0.1705 & 0.1705 & 0.1693 & 0.1692\\
        \hline
        %$\beta_2$ & 0 & 0 & 0 & 0 & 0 & 0 & 5 & 53 & 11\\
        $\beta_1$ & 0 & 0     & 0     & 24   & 841 & 841 & 1321 & 1801 & 2761\\
        $\beta_0$ & 0 & 120 & 600 & 144 & 1     & 481  & 481   & 1       & 1\\
        \end{tabular}
\caption{(Color online.) Betti numbers for the first several changes in topology.}
\label{fig:Betti}
\end{figure}

The index-$0$ critical points correspond to vertices of this graph, and the index-$1$ critical points correspond to edges.  Since we have already computed the number of connected components $\beta_0$ in the dendrogram above, the Euler relation 
$$v - e = \beta_0 - \beta_1.$$
allows us to compute $\beta_1$ --- the results in Figure \ref{fig:Betti}.  While the topology changes from the two critical points with the same radius $r \approx 0.1705$ are simultaneous, their indices are distinct and their contribution to the Betti numbers are computed separately.  The calculation does not continue for $r < 0.1686$ since we only have bounds on the Betti numbers after this point --- with the attachment of $2$-cells, the Euler relation $$v-e+f =\beta_0 - \beta_1 + \beta_2$$ alone is not sufficient to find $\beta_1$.

From Figure \ref{fig:Betti}, we see that for $0.1686 < r < 0.1692$ the configuration space is homotopy equivalent to a bouquet of $2761$ circles.  Since we know that no more $1$-cells are attached for $r < 0.1692$, we have that $\beta_1$ is monotone decreasing as the radius decreases from this point on.  Meanwhile, for small enough radius, $\beta_1 = 10$ since this is the Betti number for the configuration space of $5$ points in the plane.  Therefore, for five disks in the unit square we know that $\beta_1$ is a unimodal function of the radius and that its peak value is $\beta_1 = 2761$.

\begin{figure}
\begin{centering}
\includegraphics[width=4in]{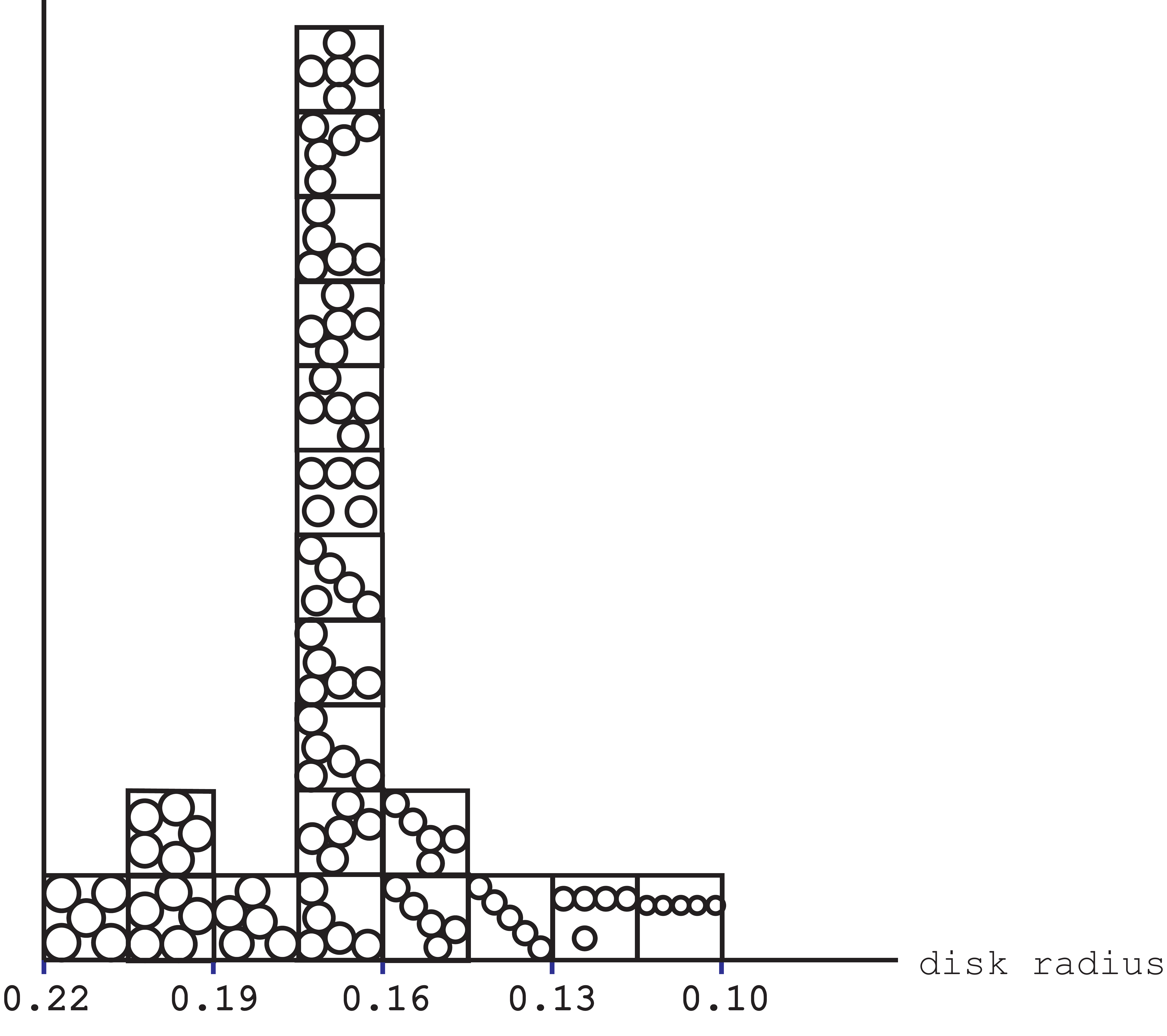}
\end{centering}
\caption{Histogram of the nineteen types of non-degenerate critical points}
\label{fig:radius5}
\end{figure}

The non-degenerate critical points are shown sorted by radius (and in rows by index) in Figure \ref{fig:radius5}. 
In this histogram, it is seen that all the critical points have radius between $0.100$ and $0.220$ --- however half of them occur between radius $0.160$ and $0.175$.  This suggests the intriguing possibility that a large number of critical points are concentrated in a small interval as the number of particles $n \to \infty$, potentially the type of situation required by the Topological Hypothesis to drive a phase transition.

We do not know with certainty that we have a complete list of all types of critical points for $5$ disks in a box --- for example, it is easily possible that some points have very small basins of attraction for the algorithms we used, and that we also failed to find them by hand.  If so, it would change the calculations we have made --- nevertheless, we expect that we have most of them and one or two missing types would likely not significantly alter the qualitative behavior we observe.  

\section{Comments}

It is probably within computational reach to implement experiments like these for greater numbers of disks, spheres in higher dimensions, and in other types of bounded regions. Flat tori and sphere are especially attractive boundaryless $2$-dimensional cases. Even if it is difficult to find a complete list of critical points for large numbers of disks, it would still be interesting to sample a large selection of them.  It might be possible, for example, to find a geometrically feasible bond network and then increase the disk radius, preserving tangency at all time, until a critical point is found.  Compiling statistics on radius and index might already give interesting information, as in Figure \ref{fig:radius5}.  New techniques for exhausting index-$0$ jammed configurations, and over all possible lattices, are described in \cite{sal-robust}, and these methods could likely be adapted to higher index points.

A promising area for future study is the asymptotic topology of $\config(n,r)$ as $n \to \infty$.  It is generally not even known how large one can make $r = r(n)$ and maintain the property that $\config(n,r)$ is connected as $n \to \infty$.  For ergodicity in a convex region, it is straightforward but important to ensure that the condition $r \le C / n$ for some region-dependent constant $C$ implies that $\config(n,r)$ is connected \cite{DLM, DLM2}.  However, spheres satisfying this condition are quite small, taking up only length rather than volume in the limit, and it may be that some configuration spaces are connected for much larger radius.

The question of threshold radius for connectivity depends on the bounding region --- the case of a circle implies that the $C/ n$ bound is the best possible when the region is not constrained.  It has been conjectured that for a flat torus though, the least dense collectively jammed configuration is the (reinforced) Kagom\'e lattice \cite{Donev}.  If an analogue of this conjecture holds for index zero points and index one critical points, for example, then $r =O(n^{-1/d})$ would ensure connectivity of $\config(n,r)$.

\begin{figure}
        \begin{tabular}{r|rrrrrrrrr}
       && \includegraphics[width=0.3in]{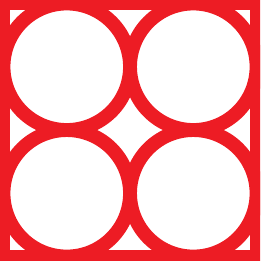}&\includegraphics[width=0.3in]{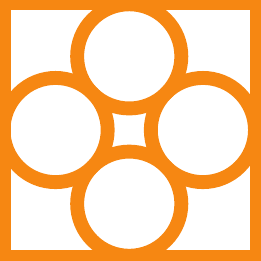}&\includegraphics[width=0.3in]{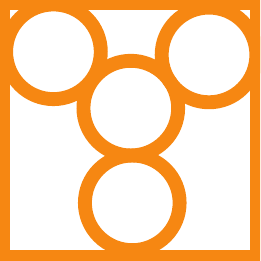}&\includegraphics[width=0.3in]{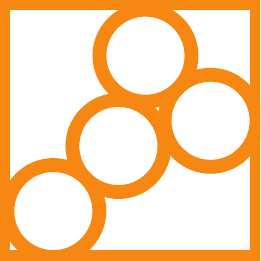} &\includegraphics[width=0.3in]{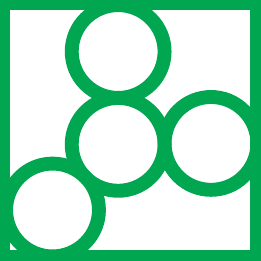}&\includegraphics[width=0.3in]{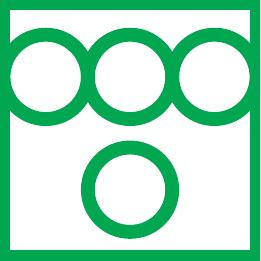}&\includegraphics[width=0.3in]{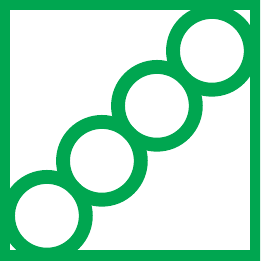}&\includegraphics[width=0.3in]{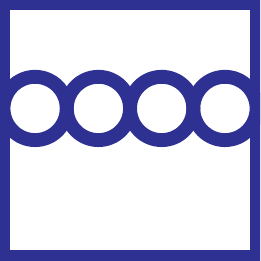}\\
       \hline
        $\beta_3$ & 0 & 0 & 0 & 0 & 0 & 0 & 0 & 0 & \fbox{6}\\
        $\beta_2$ & 0 & 0 & 0 & 0 & 0 & 0 & 5 & \fbox{53} & 11\\
        $\beta_1$ & 0 & 0 & 6 & 97 & \fbox{193} & 97 & 6 & 6 & 6\\
        $\beta_0$ & 0 & \fbox{24} & 6 & 1 & 1 & 1 & 1 & 1 & 1\\
        \end{tabular}
        \caption{(Color online.) For four disks in a square, the topology changes eight times as the radius varies.  Each row and column is unimodal, the maximum of each row is boxed.}
        \label{fig:Betti4}
\end{figure}

A more general question is the asymptotic behavior of the Betti numbers when the radius is set to $r=r(n)$ and $n \to \infty$.  We can guess but not yet prove several of the qualitative features of the asymptotic topology of $\config(n,r)$ for this situation.  For example, we expect $\beta_0$ to be roughly decreasing with radius, and each Betti number $\beta_k$ with $k \ge 1$ to be roughly unimodal in the disk radius.  In support of this conjecture, we display the Betti numbers we computed for four disks in a square in Figure \ref{fig:Betti4} --- in this case there are fewer changes in the topology, the index of a critical point is monotone decreasing in the radius, and we were able to compute all the Betti numbers.

\section*{{\bf Acknowledgements}}

We thank Yuliy Baryshnikov, Fred Cohen, Bob MacPherson, and Sal Torquato for helpful and inspiring conversations, and especially Persi Diaconis, who suggested studying the topology of these configuration spaces.  JM and MK thank the Institute for Advanced Study for supporting them in 2010--2011 when much of this work was completed.  Stanford supported Jackson Gorham as an undergraduate researcher in 2008. 

\bibliographystyle{unsrt}
\bibliography{hier-rev}

\end{document}